\documentclass[12pt,a4paper]{article}
\usepackage{amssymb,amsmath}

\usepackage[top=2.0cm,bottom=2.0cm,left=2.0cm,right=2.0cm]{geometry}
\usepackage{hyperref}
\usepackage[numbers,sort]{natbib}
\usepackage[ngerman,english]{babel}

\begin{document}

\begin{center}
{\LARGE\bf
Just another method for generating 
}
\vskip 0.25cm
{\LARGE\bf
series of even powers of number Pi
}
$$
$$
\large
{Alois Schiessl}
\centerline{}
{\tt E-Mail: aloisschiessl@web.de}
\selectlanguage{english}
\end{center}
\centerline{}
\begin{abstract}
In our previous publication we have shown a method for calculating series of even powers of $\pi$ based on the product representation of the $sinc$ function. We refer the readers to \cite{AS} for more details. In this work we apply the method to the product representation of the $cosine$ function and and thereby derive nice series formulas for even powers of the number $\pi$, such as
\[
\frac{1}{2\,!}\left(\frac{\pi}{2}\right)^2
=\sum_{\ell_1=1}^{\infty}\frac{1}{\left(2\,\ell_1-1\right)^2}
\;;\quad
\frac{1}{4\,!}\left(\frac{\pi}{2}\right)^4
=\sum_{\ell_{{2}}=2}^{\infty}\left(\sum _{\ell_{{1}}=1}^{\ell_{{2}}-1}\frac{1}{
\left(2\,\ell_1-1\right)^2\cdot
\left(2\,\ell_2-1\right)^2}
\right)\;;
\]
\[
\frac{1}{6\,!}\left(\frac{\pi}{2}\right)^6
=\sum_{\ell_{{3}}=3}^{\infty}
\left(\sum_{\ell_{{2}}=2}^{\ell_{{3}}-1}
\left(\sum_{\ell_{{1}}=1}^{\ell_{{2}}-1}
\frac{1}{
\left(2\,\ell_1-1\right)^2\cdot
\left(2\,\ell_2-1\right)^2\cdot
\left(2\,\ell_3-1\right)^2}
\right)\right)
\]
$ \\ $
Many of these formulas do not seem to be widely known.
$ \\ \\ $
$ \\ \\ $
\selectlanguage{ngerman}
\qquad\qquad\qquad\qquad\qquad\qquad\qquad\quad\textbf{Zusammenfassung}
$ \\ \\ $
In unserer früheren Publikationen haben wir ein Verfahren vorgestellt, das die Berechnung von Reihen für geradzahlige $\pi$-Potenzen unter Verwendung der $sinc$-Funktion ermöglicht. Wir verweisen die versierte Leserschaft auf \cite{AS} für nähere Details. In dieser Abhandlung wenden wir das Verfahren auf die Produktdarstellung der $cosinus$-Funktion an und erhalten weitere Reihendarstellungen für geradzahlige $\pi$-Potenzen. Die meisten der vorgestellten Reihen scheinen nicht so bekannt zu sein.
\end{abstract}
\centerline{Deutsche Version ab Seite 10}
\selectlanguage{english}
\centerline{}
\centerline{}
\centerline{***********************}
\centerline{\it In celebration of Pi Day}
\centerline{\it 3 - 14 - 2025}
\centerline{***********************}

\section*{1  Introduction}
We will start with a brief summary of the main results in \cite{AS}. Then we use the same method to establish further series representations of $\pi$-powers with even exponents. 
$ \\ \\ $
Let $M\in\mathbb N$ and $x_{1},\ldots ,x_{M}$ be $M$ independent indeterminate. We consider the product
\begin{align*}
\left(1+x_1 t\right)\cdots\left(1+x_M t\right)=\prod_{k=1}^{M}(1+x_{k}t)
\end{align*}
By expanding it, collecting and rearranging the terms, we have the identity
\begin{align*}
\prod_{k=1}^{M}(1+x_{k}t)
=1+\sigma_{M,1}t+\sigma_{M,2}t^2+\sigma_{M,3}t^3+\ldots+\sigma_{M,M}t^M=1+\sum_{k=1}^{M}\sigma_{M,k}\cdot t^{k}
\end{align*}
with the elementary symmetric polynomials
\begin{align*}
\sigma_{M,k}&=\begin{cases}
{\sum\limits_{1 \le {\ell_1} <\cdots < {\ell_k} \le M} {x_{{\ell_1}}\cdots x_{{\ell_k}}} } & (1\leq k\leq M) \\
\\
\qquad 0 & (M<k)
\end{cases}
\end{align*}
respectively.
Finally we represent the elementary symmetric polynomials by nested sums:
\[
\sigma_{M,k}
=\underbrace {\sum\limits_{{\ell_k} = k}^M {\left( {\sum\limits_{{\ell_{k - 1}} = k - 1}^{{\ell_k} - 1}  \cdots  \left( {\sum\limits_{{\ell_2} = 2}^{{\ell_3} - 1} {\left( {\sum\limits_{{\ell_1} = 1}^{{\ell_{2} - 1}} x_{\ell_1}\cdot x_{\ell_2} \cdots x_{\ell_{k-1}} \cdot x_{\ell_k} } \right)} } \right) \cdots } \right)} }_{k\;nested\;sums}
\]
Thus we can write
\[\prod_{k=1}^{M}(1+x_{k}t)=1
+\left[\sum_{\ell_1=1}^{M}x{_{\ell_1}}\right] t
\]
\[
+\left[\sum _{\ell_{{2}}=2}^{M} \left( \sum _{\ell_{{1}}=1}^{\ell_{{2}}-1}x_{{\ell_{{1}}}
} x_{{\ell_{{2}}}} \right)\right] t^2
\]
\[
+\left[\sum _{\ell_{{3}}=3}^{M} \left( \sum _{\ell_{{2}}=2}^{\ell_{{3}}-1} \left( 
\sum _{\ell_{{1}}=1}^{\ell_{{2}}-1}x_{{\ell_{{1}}}} x_{{\ell_{{2}}}} x_{{\ell_{{3}}}} \right)  \right)\right] t^3
\]
\[
+\left[{\sum _{\ell_{{4}}=4}^{M} \left(\sum _{\ell_{{3}}=3}^{\ell_4-1} \left( \sum _{\ell_{{2}}=2}^{\ell_{{3}}-1} \left( 
\sum _{\ell_{{1}}=1}^{\ell_{{2}}-1} x_{{\ell_{{1}}}}  x_{{\ell_{{2}}}} x_{{\ell_3}} x_{{\ell_4}} \right)  \right)\right)}\right] t^4\,+\ldots
\]
\[
+\left[\sum\limits_{{\ell_M} = M}^M {\left( {\sum\limits_{{\ell_{M - 1}} = M - 1}^{{\ell_M} - 1}  \cdots \left( {\sum\limits_{{\ell_2} = 2}^{{\ell_3} - 1} {\left( {\sum\limits_{{\ell_1} = 1}^{{\ell_{2} - 1}} x_{\ell_1} x_{\ell_2} \cdots x_{\ell_{M-1}} x_{\ell_M} } \right)} } \right) \cdots } \right)}\right]t^M
\]
Suitable substitutions on the left hand side and on the right hand side result in the series of $\pi$- powers presented in \cite{AS}.

\section*{2  From polynomial to product representation of cosine}
The aim of this paper is the derivation of new $\pi$-series. We consider the equation from above and substitute on the left hand side 
\[
t=-x^2\,; \quad x_k=\frac{1}{(k-\frac{1}{2})^2}\,;k\in\mathbb \{1,\ldots,n\}
\]
On the left hand side we get
\[
\prod_{k=1}^{M}\left(1-\frac{x^2}{\left(k-\frac{1}{2}\right)^2}\right)
\]
$ \\ $
That means on the right hand we have to substitute
\[
t=-x^2\,;\quad x_{\ell_k}=\frac{1}{\left(\ell_k-\frac{1}{2}\right)^2}\,;k\in\mathbb \{1,\ldots,n\}
\]
That results in
\[
1-\left[\sum_{\ell_1=1}^{M}
\frac{1}{\left(\ell_1-\frac{1}{2}\right)^2}
\right] x^2
\]
\[+\left[\sum _{\ell_{{2}}=2}^{M}\left( \sum _{\ell_{{1}}=1}^{\ell_{{2}}-1}
\frac{1}{\left(\ell_1-\frac{1}{2}\right)^2\cdot
\left(\ell_2-\frac{1}{2}\right)^2}
\right)\right] x^4
\]
\[
-\left[\sum _{\ell_{{3}}=3}^{M} \left( \sum _{\ell_{{2}}=2}^{\ell_{{3}}-1} \left( 
\sum _{\ell_{{1}}=1}^{\ell_{{2}}-1}
\frac{1}
{
\left(\ell_1-\frac{1}{2}\right)^2\cdot
\left(\ell_2-\frac{1}{2}\right)^2\cdot
\left(\ell_3-\frac{1}{2}\right)^2}
\right)\right)\right] x^6;
\]
\[
+\left[\sum _{\ell_{{4}}=4}^{M} \left( \sum _{\ell_{{3}}=3}^{\ell_{{4}}-1}\left( \sum _{\ell_{{2}}=2}^{\ell_{{3}}-1} \left( \sum _{\ell_{{1}}=1}^{\ell_{{2}
}-1}
\frac{1}
{
\left(\ell_1-\frac{1}{2}\right)^2\cdot
\cdots
\left(\ell_4-\frac{1}{2}\right)^2
}
\right)\right)\right)\right]x^8\quad -+\cdots
\]
\[
+\left(-1\right)^M
\left[\sum\limits_{{\ell_M} = M}^M {\left( {\sum\limits_{{\ell_{M - 1}} = M - 1}^{{\ell_M} - 1}  \cdots  \left( {\sum\limits_{{\ell_2} = 2}^{{\ell_3} - 1} {\left( {\sum\limits_{{\ell_1} = 1}^{{\ell_{2 - 1}}} \frac{1}
{
\left(\ell_1-\frac{1}{2}\right)^2\cdot
\cdots
\left(\ell_M-\frac{1}{2}\right)^2
} } \right)} } \right) \cdots } \right)}\right]x^{2M}
\]
$ \\ $
The fractions
\[
\frac{1}{\left(\ell_k-\frac{1}{2}\right)^2}\;;\;k\in\mathbb \{1,\ldots,M\}
\]
can be simplified as
\[
\frac{1}{\left(\ell_k-\frac{1}{2}\right)^2}
=\frac{2^2}{\left(2\,\ell_k-1\right)^2}\;;\;k\in\mathbb \{1,\ldots,M\}
\]
$ \\ $
Finally we obtain the following representation
\[
1-\left[\sum_{\ell_1=1}^{M}
\frac{2^2}{\left(2\,\ell_1-1\right)^2}
\right] x^2
\]
\[
+\left[\sum _{\ell_{{2}}=2}^{M}\left( \sum _{\ell_{{1}}=1}^{\ell_{{2}}-1}
\frac{2^2}{\left(2\,\ell_1-1\right)^2}
\cdot
\frac{2^2}{\left(2\,\ell_2-1\right)^2}
\right)\right] x^4
\]
\[
-\left[\sum _{\ell_{{3}}=3}^{M}
\left(\sum_{\ell_{{2}}=2}^{\ell_{{3}}-1} \left( 
\sum _{\ell_{{1}}=1}^{\ell_{{2}}-1}
\frac{2^2}{\left(2\,\ell_1-1\right)^2}
\cdot
\frac{2^2}{\left(2\,\ell_2-1\right)^2}
\cdot
\frac{2^2}{\left(2\,\ell_3-1\right)^2}
\right)\right)\right] x^6\;
\]
\[
+\left[\sum_{\ell_{{4}}=4}^{M}
\left(\sum_{\ell_{{3}}=3}^{\ell_{{4}}-1}
\left(\sum_{\ell_{{2}}=2}^{\ell_{{3}}-1}
\left(\sum_{\ell_{{1}}=1}^{\ell_{{2}}-1}
\frac{2^2}{\left(2\,\ell_1-1\right)^2}
\cdot
\frac{2^2}{\left(2\,\ell_2-1\right)^2}
\cdot
\frac{2^2}{\left(2\,\ell_3-1\right)^2}
\cdot
\frac{2^2}{\left(2\,\ell_4-1\right)^2}
\right)\right)\right)\right]x^8-+\cdots
\]
\[
+\left(-1\right)^M
\left[\sum\limits_{{\ell_M}=M}^M{
\left({\sum\limits_{{\ell_{M-1}}=M-1}^{{\ell_M}-1}\cdots
\left({\sum\limits_{{\ell_2}=2}^{{\ell_3}-1}{
\left({\sum\limits_{{\ell_1}=1}^{{\ell_{2-1}}}
\frac{2^2}{\left(2\,\ell_1-1\right)^2}
\cdots
\frac{2^2}{\left(2\,\ell_{M}-1\right)^2}
}\right)}}\right)\cdots}\right)}\right]x^{2\,M}
\]
$ \\ \\ $
The powers
\[2^2\,;\quad 2^2\cdot 2^2=2^4\,;\quad 2^2\cdot 2^2\cdot 2^2=2^6\,;\quad\ldots\quad\,;\,\underbrace{2^2\cdot\ldots\cdot 2^2}_{M\times\, 2^2}=2^{2\,M}
\]
in the sums don't depend on $M$ and can be excluded. This results in
\[
1-2\,^2\cdot\left[\sum_{\ell_1=1}^{M}
\frac{1}{\left(2\,\ell_1-1\right)^2}
\right] x^2
\]
\[
+2\,^4
\left[\sum _{\ell_{{2}}=2}^{M}\left( \sum _{\ell_{{1}}=1}^{\ell_{{2}}-1}
\frac{1}{\left(2\,\ell_1-1\right)^2\cdot\left(2\,\ell_2-1\right)^2}
\right)\right] x^4
\]
\[
-2\,^6
\left[\sum _{\ell_{{3}}=3}^{M}
\left(\sum_{\ell_{{2}}=2}^{\ell_{{3}}-1} \left( 
\sum _{\ell_{{1}}=1}^{\ell_{{2}}-1}
\frac{1}{\left(2\,\ell_1-1\right)^2
\cdot
\left(2\,\ell_2-1\right)^2
\cdot
\left(2\,\ell_3-1\right)^2}
\right)\right)\right] x^6\;
\]
\[
+2\,^8
\left[\sum_{\ell_{{4}}=4}^{M}
\left(\sum_{\ell_{{3}}=3}^{\ell_{{4}}-1}
\left(\sum_{\ell_{{2}}=2}^{\ell_{{3}}-1}
\left(\sum_{\ell_{{1}}=1}^{\ell_{{2}}-1}
\frac{1}{\left(2\,\ell_1-1\right)^2}
\cdot
\frac{1}{\left(2\,\ell_2-1\right)^2
\cdot
\left(2\,\ell_3-1\right)^2
\cdot
\left(2\,\ell_4-1\right)^2}
\right)\right)\right)\right]x^8-+\cdots
\]
\[
+\left(-1\right)^M\cdot2^{\,2\,M}
\left[\sum\limits_{{\ell_M}=M}^M{
\left({\sum\limits_{{\ell_{M-1}}=M-1}^{{\ell_M}-1}\cdots
\left({\sum\limits_{{\ell_2}=2}^{{\ell_3}-1}{
\left({\sum\limits_{{\ell_1}=1}^{{\ell_{2-1}}}
\frac{1}{\left(2\,\ell_1-1\right)^2
\cdots
\left(2\,\ell_{M}-1\right)^2}
}\right)}}\right)\cdots}\right)}\right]x^{2\,M}
\]
$ \\ $
Now let us examine when taking limit $M \to \infty \,$. The left hand side results in an infinite product
\[
\underset{n\to \infty}{\mathop{\lim }}\;\;
\prod_{k=1}^{M}\left(1-\frac{x^2}{\left(k-\frac{1}{2}\right)^2}\right)
=\prod_{k=1}^{\infty}\left(1-\frac{x^2}{\left(k-\frac{1}{2}\right)^2}\right),
\]
$ \\ $
which we identify as the product representation of the cosine function
\[
\cos\left(\pi x\right)=\prod_{k=1}^{\infty}\left(1-\frac{x^2}{\left(k-\frac{1}{2}\right)^2}\right)
\]
The power series of $\cos\left(\pi\,x\right)$ can be found in any formula collection:
\[
\cos\left(\pi\,x\right)
=1-\frac{\pi^2}{2\,!}x^2+\frac{\pi^4}{4\,!}x^4-\frac{\pi^6}{6\,!}x^6+\frac{\pi^8}{8\,!}x^8-+\cdots
\]
On the right hand side when taking limit $M \to \infty $, we obtain an infinite power series with coefficients consisting of infinite nested sums
\[
1-2^2\cdot\left[\sum_{\ell_1=1}^{\infty}
\frac{1}{\left(2\,\ell_1-1\right)^2}
\right] x^2
\]
\[
+2\,^4
\left[\sum _{\ell_{{2}}=2}^{\infty}\left( \sum _{\ell_{{1}}=1}^{\ell_{{2}}-1}
\frac{1}{\left(2\,\ell_1-1\right)^2\cdot\left(2\,\ell_2-1\right)^2}
\right)\right] x^4
\]
\[
-2\,^6
\left[\sum _{\ell_{{3}}=3}^{\infty}
\left(\sum_{\ell_{{2}}=2}^{\ell_{{3}}-1} \left( 
\sum _{\ell_{{1}}=1}^{\ell_{{2}}-1}
\frac{1}{\left(2\,\ell_1-1\right)^2
\cdot
\left(2\,\ell_2-1\right)^2
\cdot
\left(2\,\ell_3-1\right)^2}
\right)\right)\right] x^6\;
\]
\[
+2\,^8
\left[\sum_{\ell_{{4}}=4}^{\infty}
\left(\sum_{\ell_{{3}}=3}^{\ell_{{4}}-1}
\left(\sum_{\ell_{{2}}=2}^{\ell_{{3}}-1}
\left(\sum_{\ell_{{1}}=1}^{\ell_{{2}}-1}
\frac{1}{\left(2\,\ell_1-1\right)^2}
\cdot
\frac{1}{\left(2\,\ell_2-1\right)^2
\cdot
\left(2\,\ell_3-1\right)^2
\cdot
\left(2\,\ell_4-1\right)^2}
\right)\right)\right)\right]x^8\;-+\cdots
\]
\section*{3  From cosine product representation to series of Pi}
We therefore have two different representations with power series for $\cos\left(\pi\,x\right)$. One is the classical power series
\[
\cos\left(\pi\,x\right)=1-\frac{\pi^2}{2!}x^2+\frac{\pi^4}{4!}x^4-\frac{\pi^6}{6!}x^6+-\cdots
\]
The second is a representation by nested sums. According to the uniqueness theorem for power series, both must be identical. The following must therefore apply:
\[
1-\frac{\pi^2}{2!}x^2+\frac{\pi^4}{4!}x^4-\frac{\pi^6}{6!}x^6
-\frac{\pi^8}{8\,!}x^8-+\cdots
\]
\[
=1-2\,^2\left[\sum_{\ell_1=1}^{\infty}
\frac{1}{\left(2\,\ell_1-1\right)^2}
\right] x^2
\]
\[
+2\,^4
\left[\sum _{\ell_{{2}}=2}^{\infty}\left( \sum _{\ell_{{1}}=1}^{\ell_{{2}}-1}
\frac{1}{\left(2\,\ell_1-1\right)^2\cdot\left(2\,\ell_2-1\right)^2}
\right)\right] x^4
\]
\[
-2\,^6
\left[\sum _{\ell_{{3}}=3}^{\infty}
\left(\sum_{\ell_{{2}}=2}^{\ell_{{3}}-1} \left( 
\sum _{\ell_{{1}}=1}^{\ell_{{2}}-1}
\frac{1}{\left(2\,\ell_1-1\right)^2
\cdot
\left(2\,\ell_2-1\right)^2
\cdot
\left(2\,\ell_3-1\right)^2}
\right)\right)\right] x^6\;
\]
\[
+2\,^8
\left[\sum_{\ell_{{4}}=4}^{\infty}
\left(\sum_{\ell_{{3}}=3}^{\ell_{{4}}-1}
\left(\sum_{\ell_{{2}}=2}^{\ell_{{3}}-1}
\left(\sum_{\ell_{{1}}=1}^{\ell_{{2}}-1}
\frac{1}{\left(2\,\ell_1-1\right)^2}
\cdot
\frac{1}{\left(2\,\ell_2-1\right)^2
\cdot
\left(2\,\ell_3-1\right)^2
\cdot
\left(2\,\ell_4-1\right)^2}
\right)\right)\right)\right]x^8\;-+\cdots
\]
$ \\ $
If two power series are equal, then the corresponding coefficients must be equal.
$ \\ \\ $
On equating the coefficients of power $x^2$ we obtain
\[
\frac{\pi^2}{2!}
=2^2\sum_{\ell_1=1}^{\infty}\frac{1}{\left(2\,\ell_1-1\right)^2}
\]
Dividing by $2^2$, we get
\[
\frac{\pi^2}{2\,!\cdot 2^2}
=\frac{1}{2\,!}\left(\frac{\pi}{2}\right)^2
=\frac{\pi^2}{8}
=\sum_{\ell_1=1}^{\infty}\frac{1}{\left(2\,\ell_1-1\right)^2}
=1+\frac{1}{3^2}+\frac{1}{5^2}+\frac{1}{7^2}+\frac{1}{9^2}+\ldots
\]
This series often occurs in connection with solutions to the Basel problem.
$ \\ \\ $
On equating the coefficients of power $x^4$, we find on the left side a rational multiple of $\pi^4$ and on the right hand side a double sum
\[
\frac{\pi^4}{4\,!}
=2^4\sum_{\ell_{{2}}=2}^{\infty}\left(\sum _{\ell_{{1}}=1}^{\ell_{{2}}-1}\frac{1}{
\left(2\,\ell_1-1\right)^2\cdot
\left(2\,\ell_2-1\right)^2}
\right)
\]
Dividing by $2^4$ gives
\[
\frac{\pi^4}{4\,!\,2^4}
=\frac{1}{4\,!}\left(\frac{\pi}{2}\right)^4
=\sum_{\ell_{{2}}=2}^{\infty}\left(\sum _{\ell_{{1}}=1}^{\ell_{{2}}-1}\frac{1}{
\left(2\,\ell_1-1\right)^2\cdot
\left(2\,\ell_2-1\right)^2}
\right)
\]
We define the right hand side limit 
\[
\sum_{\ell_{{2}}=2}^{\infty}\left(\sum _{\ell_{{1}}=1}^{\ell_{{2}}-1}\frac{1}{
\left(2\,\ell_1-1\right)^2\cdot
\left(2\,\ell_2-1\right)^2}
\right)
\]
as the limit of the partial sums
\[
\sum_{\ell_2=2}^{2}
\left(\sum_{\ell_1=1}^{\ell_2-1}
\frac{1}
{
\left(2\ell_1-1\right)^2
\left(2\ell_2-1\right)^2
}
\right)
=\frac{1}{1^2\cdot 3^2}
\]
\[
\sum_{\ell_2=2}^{3}
\left(\sum_{\ell_1=1}^{\ell_2-1}
\frac{1}
{
\left(2\ell_1-1\right)^2
\left(2\ell_2-1\right)^2
}
\right)
=\frac{1}{1^2\cdot 3^2}+\frac{1}{1^2\cdot 5^2}+\frac{1}{3^2\cdot 5^2}
\]
\[
\sum_{\ell_2=2}^{4}
\left(\sum_{\ell_1=1}^{\ell_2-1}
\frac{1}
{
\left(2\ell_1-1\right)^2
\left(2\ell_2-1\right)^2
}
\right)
=\frac{1}{1^2\cdot 3^2}+\frac{1}{1^2\cdot 5^2}+\frac{1}{1^2\cdot7^2}
+\frac{1}{3^2\cdot 5^2}+\frac{1}{3^2\cdot7^2}
+\frac{1}{5^2\cdot7^2}
\]
\centerline{\vdots}
$ $
Next, equating the coefficients of power $x^6$ gives a rational multiple of $\pi^6$ and a triple sum:
\[
\frac{\pi^6}{6\,!}
=2^6\sum_{\ell_{{3}}=3}^{\infty}
\left(\sum_{\ell_{{2}}=2}^{\ell_{{3}}-1}
\left(\sum_{\ell_{{1}}=1}^{\ell_{{2}}-1}
\frac{1}{
\left(2\,\ell_1-1\right)^2\cdot
\left(2\,\ell_2-1\right)^2\cdot
\left(2\,\ell_3-1\right)^2}
\right)\right)
\]
We can also simplify this into
\[
\frac{\pi^6}{6\,!\cdot 2^6}
=\frac{1}{6\,!}\left(\frac{\pi}{2}\right)^6
=\sum_{\ell_{{3}}=3}^{\infty}
\left(\sum_{\ell_{{2}}=2}^{\ell_{{3}}-1}
\left(\sum_{\ell_{{1}}=1}^{\ell_{{2}}-1}
\frac{1}{
\left(2\,\ell_1-1\right)^2\cdot
\left(2\,\ell_2-1\right)^2\cdot
\left(2\,\ell_3-1\right)^2}
\right)\right)
\]
\\
We also define the limit
\[
\sum_{\ell_{{3}}=3}^{\infty}
\left(\sum_{\ell_{{2}}=2}^{\ell_{{3}}-1}
\left(\sum_{\ell_{{1}}=1}^{\ell_{{2}}-1}
\frac{1}{
\left(2\,\ell_1-1\right)^2\cdot
\left(2\,\ell_2-1\right)^2\cdot
\left(2\,\ell_3-1\right)^2}
\right)\right)
\]
as the limit of the partial sums
\[
\sum_{\ell_3=3}^{3}
\left(\sum_{\ell_2=2}^{\ell_3-1}
\left(\sum_{\ell_1=1}^{\ell_2-1}
\frac{1}
{
\left(2\ell_1-1\right)^2
\left(2\ell_2-1\right)^2
\left(2\ell_3-1\right)^2
}
\right)\right)
=\frac{1}{1^{2} \cdot 3^{2}\cdot 5^{2}}
\]
\centerline{}
\[
\sum_{\ell_3=3}^{4}
\left(\sum_{\ell_2=2}^{\ell_3-1}
\left(\sum_{\ell_1=1}^{\ell_2-1}
\frac{1}
{
\left(2\ell_1-1\right)^2
\left(2\ell_2-1\right)^2
\left(2\ell_3-1\right)^2
}
\right)\right)
\]
\[
=\frac{1}{1^{2} \cdot 3^{2} \cdot 5^{2}}+\frac{1}{1^{2} \cdot 3^{2} \cdot 7^{2}}+\frac{1}{1^{2} \cdot 5^{2} \cdot 7^{2}}+\frac{1}{3^{2} \cdot 5^{2} \cdot 7^{2}}
\]
\centerline{}
\[
\sum_{\ell_3=3}^{5}
\left(\sum_{\ell_2=2}^{\ell_3-1}
\left(\sum_{\ell_1=1}^{\ell_2-1}
\frac{1}
{
\left(2\ell_1-1\right)^2
\left(2\ell_2-1\right)^2
\left(2\ell_3-1\right)^2
}
\right)\right)
\]
\[
=\frac{1}{1^{2} \cdot 3^{2} \cdot 5^{2}}+\frac{1}{1^{2} \cdot 3^{2} \cdot 7^{2}}+\frac{1}{1^{2} \cdot 3^{2} \cdot 9^{2}}+\frac{1}{1^{2} \cdot 5^{2} \cdot 7^{2}}+\frac{1}{1^{2} \cdot 5^{2} \cdot 9^{2}}
\]
\[
+\frac{1}{1^{2} \cdot 7^{2} \cdot 9^{2}}+\frac{1}{3^{2} \cdot 5^{2} \cdot 7^{2}}+\frac{1}{3^{2} \cdot 5^{2} \cdot 9^{2}}+\frac{1}{3^{2} \cdot 7^{2} \cdot 9^{2}}+\frac{1}{5^{2} \cdot 7^{2} \cdot 9^{2}}
\]
$ \\ \\ $
Continuing with equating coefficients of powers $x^8$ and $x^{10}$ yields
\[
\frac{1}{8!}\left(\frac{\pi}{2}\right)^8=
\]
\[=
\sum_{\ell_4=4}^{\infty}
\left(\sum_{\ell_3=3}^{\ell_4-1}
\left(\sum_{\ell_2=2}^{\ell_3-1}
\left(\sum_{\ell_1=1}^{\ell_2-1}
\frac{1}
{
\left(2\ell_1-1\right)^2
\left(2\ell_2-1\right)^2
\left(2\ell_3-1\right)^2
\left(2\ell_4-1\right)^2
}
\right)\right)\right)
\]
\centerline{}
\[
\frac{1}{10!}\left(\frac{\pi}{2}\right)^{10}=
\]
\[
=\sum_{\ell_5=5}^{\infty}
\left(\sum_{\ell_4=4}^{\ell_5-1}
\left(\sum_{\ell_3=3}^{\ell_4-1}
\left(\sum_{\ell_2=2}^{\ell_3-1}
\left(\sum_{\ell_1=1}^{\ell_2-1}
\frac{1}
{
\left(2\ell_1-1\right)^2
\left(2\ell_2-1\right)^2
\left(2\ell_3-1\right)^2
\left(2\ell_4-1\right)^2
\left(2\ell_5-1\right)^2
}
\right)\right)\right)\right)
\]
$ \\ \\ $
Now let's examine the  general case equating $x^{\,2\,n}$ with $n\in\mathbb N$. Thus we have
\[
\frac{\pi^{2\,n}}{\left(2\,n\right)!}
\]
\[
=2\,^{2\,n}\sum\limits_{{\ell_n}= n}^\infty {\left( {\sum\limits_{{\ell_{n-1}}=n-1}^{{\ell_n}-1} \cdots \left({\sum\limits_{{\ell_2}=2}^{{\ell_3}-1} {\left( {\sum\limits_{{\ell_1}=1}^{{\ell_{2-1}}}\frac{1}
{
\left(2\ell_1-1\right)^2
\left(2\ell_2-1\right)^2
\cdots
\left(2\ell_{n-1}-1\right)^2
\left(2\ell_n-1\right)^2
}}\right)}}\right)\cdots}\right)}
\]
$ \\ \\ $
Next dividing the equation by $2\,^{2\,n}$ and simplifying the left hand side finally leads to the result
\[
\frac{1}{\left(2\,n\right)!}\left(\frac{\pi}{2}\right)^{2\,n}=
\]
\[=\sum\limits_{{\ell_n}= n}^\infty {\left( {\sum\limits_{{\ell_{n-1}}=n-1}^{{\ell_n}-1} \cdots \left({\sum\limits_{{\ell_2}=2}^{{\ell_3}-1} {\left( {\sum\limits_{{\ell_1}=1}^{{\ell_{2-1}}}\frac{1}
{
\left(2\ell_1-1\right)^2
\left(2\ell_2-1\right)^2
\cdots
\left(2\ell_{n-1}-1\right)^2
\left(2\ell_n-1\right)^2
}}\right)}}\right)\cdots}\right)}
\]
To conclude this section let's focus on the convergence of the sums. Since the sums only consiting of positive terms, it is easily verified by comparison test.
$ \\ \\ $
The convergence of the sum
\[
\sum_{\ell_1=1}^{\infty}\frac{1}{\left(2\,\ell_1-1\right)^2}
=1+\frac{1}{3^2}+\frac{1}{5^2}+\frac{1}{7^2}+\frac{1}{9^2}+\ldots
\]
is well known. There is nothing more to do.
$ \\ \\ $
For higher sums we use the well known sequence of series \cite{KK} \cite{RC}
\[
\sum _{\ell_{{2}}=1}^{\infty} \left( \sum _{\ell_{{1}}=1}^{\infty} 
\frac{1}
{
{\left(2\ell_1-1\right)}^2\cdot
{\left(2\ell_2-1\right)}^2
}\right)
=\left(\frac{\pi^2}{8}\right)^2;
\]
\[
\sum _{\ell_{{3}}=1}^{\infty} \left( \sum _{\ell_{{2}}=1}^{\infty} \left( \sum _
{\ell_{{1}}=1}^{\infty}
\frac{1}{
{\left(2\ell_1-1\right)}^2\cdot
{\left(2\ell_2-1\right)}^2\cdot
{\left(2\ell_3-1\right)}^2
}\right)\right)
=\left(\frac{\pi^2}{8}\right)^3;
\]
\[
\sum_{\ell_{{4}}=1}^{\infty}
\left(\sum_{\ell_{{3}}=1}^{\infty}
\left(\sum_{\ell_{{2}}=1}^{\infty}
\left(\sum_{\ell_{{1}}=1}^{\infty}
\frac{1}{
\left(2\ell_1-1\right)^2\cdot
\left(2\ell_2-1\right)^2\cdot
\left(2\ell_3-1\right)^2\cdot
\left(2\ell_4-1\right)^2}
\right)\right)\right)
=\left(\frac{\pi^2}{8}\right)^4;
\]
\[
\vdots
\]
\[
\sum\limits_{{\ell_n} = 1}^{\infty} {\left( {\sum\limits_{{\ell_{n - 1}} =1}^{\infty}  \cdots  \left( {\sum\limits_{{\ell_2} = 1}^{\infty} {\left( {\sum\limits_{{\ell_1} = 1}^{\infty} {\frac{1}{
{\left(2\ell_1-1\right)^2 \cdots
\left(2\ell_n-1\right)^2}}} } \right)} } \right) \cdots } \right)}
=\left(\frac{\pi^2}{8}\right)^n\]
$ \\ $
Comparing it with the derived series, we observe
\[
\sum _{\ell_{{2}}=2}^{\infty} \left( \sum _{\ell_{{1}}=1}^{\ell_{{2}}-1} 
\frac{1}{\left(2\ell_1-1\right)^2 \cdot
\left(2\ell_2-1\right)^2}\right)
<\sum _{\ell_{{2}}=1}^{\infty} \left( \sum _{\ell_{{1}}=1}^{\infty} 
\frac{1}{\left(2\ell_1-1\right)^2 \cdot
\left(2\ell_2-1\right)^2}\right)=\left(\frac{\pi^2}{8}\right)^2;
\]
\begin{align*}
&\sum _{\ell_{{3}}=3}^{\infty} \left( \sum _{\ell_{{2}}=2}^{\ell_{{3}}-1} \left( \sum _
{\ell_{{1}}=1}^{\ell_{{2}-1}} \frac{1}{
\left(2\ell_1-1\right)^2 \cdot
\left(2\ell_2-1\right)^2 \cdot
\left(2\ell_3-1\right)^2}\right)\right)\\ 
<
&\sum _{\ell_{{3}}=1}^{\infty} \left( \sum _{\ell_{{2}}=1}^{\infty} \left( \sum _
{\ell_{{1}}=1}^{\infty} \frac{1}{
\left(2\ell_1-1\right)^2\cdot
\left(2\ell_2-1\right)^2\cdot
\left(2\ell_3-1\right)^2}
\right)\right)=\left(\frac{\pi^2}{8}\right)^3;
\end{align*}
\centerline{}
\begin{align*}
&\sum_{\ell_{{4}}=4}^{\infty}
\left(\sum_{\ell_{{3}}=2}^{\ell_4-1}
\left(\sum_{\ell_{{2}}=2}^{\ell_3-1}
\left(\sum_{\ell_{{1}}=1}^{\ell_2-1}
\frac{1}{
\left(2\ell_1-1\right)^2\cdot
\left(2\ell_2-1\right)^2\cdot
\left(2\ell_3-1\right)^2\cdot
\left(2\ell_4-1\right)^2}
\right)\right)\right)
\\
<
&\sum_{\ell_{{4}}=1}^{\infty}
\left(\sum_{\ell_{{3}}=1}^{\infty}
\left(\sum_{\ell_{{2}}=1}^{\infty}
\left(\sum_{\ell_{{1}}=1}^{\infty}
\frac{1}{
\left(2\ell_1-1\right)^2\cdot
\left(2\ell_2-1\right)^2\cdot
\left(2\ell_3-1\right)^2\cdot
\left(2\ell_4-1\right)^2}
\right)\right)\right)
=\left(\frac{\pi^2}{8}\right)^4;
\end{align*}
\[
\vdots
\]
More generally, for $n\in\mathbb N$, we obtain
\begin{align*}
&\sum\limits_{{\ell_n} = n}^{\infty} {\left( {\sum\limits_{{\ell_{n - 1}} = n - 1}^{{\ell_n} - 1}  \cdots  \left( {\sum\limits_{{\ell_2} = 2}^{{\ell_3} - 1} {\left( {\sum\limits_{{\ell_1} = 1}^{{\ell_{2 - 1}}} {\frac{1}{{
\left(2\ell_1-1\right)^2\cdots
\left(2\ell_n-1\right)^2}}}}\right)}}\right)\cdots}\right)}
\\
<
&\sum\limits_{{\ell_n} = 1}^{\infty} {\left( {\sum\limits_{{\ell_{n - 1}} =1}^{\infty}  \cdots  \left( {\sum\limits_{{\ell_2} = 1}^{\infty} {\left( {\sum\limits_{{\ell_1} = 1}^{\infty} {\frac{1}{{
\left(2\ell_1-1\right)^2\cdots
\left(2\ell_n-1\right)^2
}}} } \right)} } \right) \cdots } \right)}
=\left(\frac{\pi^2}{8}\right)^n
\end{align*}
We conclude, that every $n^{th}$-nested sum is bounded by $\left(\frac{\pi^2}{8}\right)^n $, which implies convergence.

\newpage
\selectlanguage{ngerman}
\begin{center}
{\LARGE\bf
Ein weitere Methode zur Berechnung 
}
\vskip 0.25cm
{\LARGE\bf
von Reihen für geradzahlige Pi-Potenzen
}
$$
$$
\large
{Alois Schiessl}
\vskip 0.25cm
{\tt E-Mail: aloisschiessl@web.de}
\selectlanguage{ngerman}\date{\today}
\end{center}
$ \\ $
$ \\ $
\centerline{\textbf{Zusammenfassung}}
$ \\ $
In einer unserer früheren Veröffentlichung haben wir ein Verfahren vorgestellt, das die Berechnung von Reihen für geradzahlige $\pi$-Potenzen ermöglicht. Wir verweisen die versierte Leserschaft auf \cite{AS} für nähere Details. Wir folgten dazu Eulers Idee mit Hilfe der Produktdarstellung der $sinc$-Funktion das Basler Problems \cite{E1}, \cite{E2} zu lösen und verallgemeinerten die Idee auf höhere Potenzen. In dieser Abhandlung wenden wir das Verfahren auf die Produktdarstellung des $cosinus$ an und erhalten weitere Reihendarstellungen für geradzahlige $\pi$-Potenzen. Bemerkenswert ist hierbei, dass in diesem Fall nur ungerade Potenzen der natürlichen Zahlen im Nenner der Summanden auftreten. Wir geben einige Beispiele:
\[
\frac{1}{2!}\left(\frac{\pi}{2}\right)^2
=\sum_{\ell_1=1}^{\infty}\frac{1}{\left(2\ell_1-1\right)^2}\,;\quad
\frac{1}{4!}\left(\frac{\pi}{2}\right)^4
=
\sum_{\ell_2=2}^{\infty}
\left(\sum_{\ell_1=1}^{\ell_2-1}
\frac{1}
{
\left(2\ell_1-1\right)^2
\left(2\ell_2-1\right)^2
}
\right)\,;
\]
\[
\frac{1}{6!}\left(\frac{\pi}{2}\right)^6
=
\sum_{\ell_3=3}^{\infty}
\left(\sum_{\ell_2=2}^{\ell_3-1}
\left(\sum_{\ell_1=1}^{\ell_2-1}
\frac{1}
{
\left(2\ell_1-1\right)^2
\left(2\ell_2-1\right)^2
\left(2\ell_3-1\right)^2
}
\right)\right)
\]
$ \\ $
Die meisten der angegebenen Reihen scheinen nicht so bekannt zu sein.
$ \\ $
\centerline{}
\centerline{*********************************************}
\centerline{\it Veröffentlicht anlässlich des Tages der Kreiszahl $\pi$}

\centerline{\it 3 - 14 - 2025}
\centerline{*********************************************}
\section*{1  Einführung}
Das Verfahren ist das gleiche wie in der bereits erwähnten Veröffentlichung \cite{AS}. Diesmal verwenden wir die Produktdarstellung von $\cos\left(\pi\,x\right)$ um weitere Reihendarstellungen für geradzahlige $\pi$-Potenzen zu erhalten. Wir geben zunächst eine kurze Zusammenfassung der verwendeten Methode.
$ \\ $
$ \\ $
Zu $M\in\mathbb N$ betrachten wir die $M$ Unbestimmten $x_{1},\ldots ,x_{M}$ und bilden das Produkt
\begin{align*}
\left(1+x_1 t\right)\cdots\left(1+x_M t\right)=\prod_{k=1}^{M}(1+x_{k}t)
\end{align*}
Wir lösen die Klammern auf und fassen gleiche Potenzen zusammen, so dass wir schließlich die folgende Darstellung erhalten:
\begin{align*}
\prod_{k=1}^{M}(1+x_{k}t)
=1+\sigma_{M,1}t+\sigma_{M,2}t^2+\sigma_{M,3}t^3+\ldots+\sigma_{M,M}t^M=1+\sum_{k=1}^{M}\sigma_{M,k}\cdot t^{k}
\end{align*}
Hierbei bestehen die Koeffizienten aus den elementar-symmetrischen Polynomen
\begin{align*}
\sigma_{M,k}&=\begin{cases}
{\sum\limits_{1 \le {\ell_1} <\cdots < {\ell_k} \le M} {x_{{\ell_1}}\cdots x_{{\ell_k}}} } & (1\leq k\leq M) 
\\
\\
\qquad 0 & (M<k)
\end{cases}
\end{align*}
Wir haben die elementar-symmetrischen Polynome
\[
{\sum\limits_{1 \le {\ell_1} <\cdots < {\ell_k} \le M} {x_{{\ell_1}}\cdots x_{{\ell_k}}} }
\]
durch $k$ ineinander geschachtelte Summen dargestellt:
\[
\sigma_{M,k}
=\underbrace {\sum\limits_{{\ell_k} = k}^M {\left( {\sum\limits_{{\ell_{k - 1}} = k - 1}^{{\ell_k} - 1}  \cdots  \left( {\sum\limits_{{\ell_2} = 2}^{{\ell_3} - 1} {\left( {\sum\limits_{{\ell_1} = 1}^{{\ell_{2} - 1}} x_{\ell_1}\cdot x_{\ell_2} \cdots x_{\ell_{k-1}} \cdot x_{\ell_k} } \right)} } \right) \cdots } \right)} }_{k\;-fache\;Summe}
\]
Somit lässt sich das Produkt
\[
\prod_{k=1}^{M}(1+x_{k}t)=\left(1+x_1 t\right)\cdots\left(1+x_M t\right)
\]
durch ineinander geschachtelten Summen angeben:
\[\prod_{k=1}^{M}(1+x_{k}t)=1
+\left[\sum_{\ell_1=1}^{M}x{_{\ell_1}}\right] t
\]
\[
+\left[\sum _{\ell_{{2}}=2}^{M} \left( \sum _{\ell_{{1}}=1}^{\ell_{{2}}-1}x_{{\ell_{{1}}}
} x_{{\ell_{{2}}}} \right)\right] t^2
\]
\[
+\left[\sum _{\ell_{{3}}=3}^{M} \left( \sum _{\ell_{{2}}=2}^{\ell_{{3}}-1} \left( 
\sum _{\ell_{{1}}=1}^{\ell_{{2}}-1}x_{{\ell_{{1}}}} x_{{\ell_{{2}}}} x_{{\ell_{{3}}}} \right)  \right)\right] t^3
\]
\[
+\left[{\sum _{\ell_{{4}}=4}^{M} \left(\sum _{\ell_{{3}}=3}^{\ell_4-1} \left( \sum _{\ell_{{2}}=2}^{\ell_{{3}}-1} \left( 
\sum _{\ell_{{1}}=1}^{\ell_{{2}}-1} x_{{\ell_{{1}}}}  x_{{\ell_{{2}}}} x_{{\ell_3}} x_{{\ell_4}} \right)  \right)\right)}\right] t^4
\]
\centerline{\vdots}
\[
+\left[\sum\limits_{{\ell_M} = M}^M {\left( {\sum\limits_{{\ell_{M - 1}} = M - 1}^{{\ell_M} - 1}  \cdots \left( {\sum\limits_{{\ell_2} = 2}^{{\ell_3} - 1} {\left( {\sum\limits_{{\ell_1} = 1}^{{\ell_{2} - 1}} x_{\ell_1} x_{\ell_2} \cdots x_{\ell_{M-1}} x_{\ell_M} } \right)} } \right) \cdots } \right)}\right]t^M
\]
$ \\ $
Durch geeignete Substitutionen linksseitig und rechtsseitig und Auswertung des Grenzübergang $M \to \infty$ ergaben sich die in \cite{AS} angegeben Reihen für geradzahlige  $\pi$-Potenzen.
\section*{2  Vom Polynom zur Produktdarstellung des cosinus}
Nun kommen wir zu den neuen Reihen für $\pi$-Potenzen. Ausgangspunkt ist die bereits erwähnte Darstellung
\[
\prod_{k=1}^{M}(1+x_{k}t)=1
+\left[\sum_{\ell_1=1}^{M}x{_{\ell_1}}\right] t
\]
\[
+\left[\sum _{\ell_{{2}}=2}^{M} \left( \sum _{\ell_{{1}}=1}^{\ell_{{2}}-1}x_{{\ell_{{1}}}
} x_{{\ell_{{2}}}} \right)\right] t^2
\]
\[
+\left[\sum _{\ell_{{3}}=3}^{M} \left( \sum _{\ell_{{2}}=2}^{\ell_{{3}}-1} \left( 
\sum _{\ell_{{1}}=1}^{\ell_{{2}}-1}x_{{\ell_{{1}}}} x_{{\ell_{{2}}}} x_{{\ell_{{3}}}} \right)  \right)\right] t^3
\]
\[
+\left[{\sum _{\ell_{{4}}=4}^{M} \left(\sum _{\ell_{{3}}=3}^{\ell_4-1} \left( \sum _{\ell_{{2}}=2}^{\ell_{{3}}-1} \left( 
\sum _{\ell_{{1}}=1}^{\ell_{{2}}-1} x_{{\ell_{{1}}}}  x_{{\ell_{{2}}}} x_{{\ell_3}} x_{{\ell_4}} \right)  \right)\right)}\right] t^4
\]
\[
+\left[\sum\limits_{{\ell_M} = M}^M {\left( {\sum\limits_{{\ell_{M - 1}} = M - 1}^{{\ell_M} - 1}  \cdots \left( {\sum\limits_{{\ell_2} = 2}^{{\ell_3} - 1} {\left( {\sum\limits_{{\ell_1} = 1}^{{\ell_{2} - 1}} x_{\ell_1} x_{\ell_2} \cdots x_{\ell_{M-1}} x_{\ell_M} } \right)} } \right) \cdots } \right)}\right]t^M
\]
$ \\ $
Wir substituieren diesmal linksseitig
\[
t=-x^2\,; \quad x_k=\frac{1}{(k-\frac{1}{2})^2}\;;\;k\in\mathbb \{1,\ldots,n\}
\]
Das ergibt dann auf der linken Seite
\[
\prod_{k=1}^{M}\left(1-\frac{x^2}{\left(k-\frac{1}{2}\right)^2}\right)
\]
$ \\ $
Auf der rechten Seite müssen wir dann wie folgt substituieren:
\[
t=-x^2\,;\quad x_{\ell_k}=\frac{1}{\left(\ell_k-\frac{1}{2}\right)^2}\;;\;k\in\mathbb \{1,\ldots,n\}
\]
Das führt auf der rechten Seite zu
\[
1-\left[\sum_{\ell_1=1}^{M}
\frac{1}{\left(\ell_1-\frac{1}{2}\right)^2}
\right] x^2
\]
\[
+\left[\sum _{\ell_{{2}}=2}^{M}\left( \sum _{\ell_{{1}}=1}^{\ell_{{2}}-1}
\frac{1}{\left(\ell_1-\frac{1}{2}\right)^2\cdot
\left(\ell_2-\frac{1}{2}\right)^2}
\right)\right] x^4
\]
\[
-\left[\sum _{\ell_{{3}}=3}^{M} \left( \sum _{\ell_{{2}}=2}^{\ell_{{3}}-1} \left( 
\sum _{\ell_{{1}}=1}^{\ell_{{2}}-1}
\frac{1}
{
\left(\ell_1-\frac{1}{2}\right)^2\cdot
\left(\ell_2-\frac{1}{2}\right)^2\cdot
\left(\ell_3-\frac{1}{2}\right)^2}
\right)\right)\right] x^6\;
\]
\[
+\left[\sum _{\ell_{{4}}=4}^{M} \left( \sum _{\ell_{{3}}=3}^{\ell_{{4}}-1}\left( \sum _{\ell_{{2}}=2}^{\ell_{{3}}-1} \left( \sum _{\ell_{{1}}=1}^{\ell_{{2}
}-1}
\frac{1}
{
\left(\ell_1-\frac{1}{2}\right)^2\cdot
\cdots
\left(\ell_4-\frac{1}{2}\right)^2
}
\right)\right)\right)\right]x^8\quad -+\cdots
\]
\[
+\left(-1\right)^M
\left[\sum\limits_{{\ell_M} = M}^M {\left( {\sum\limits_{{\ell_{M - 1}} = M - 1}^{{\ell_M} - 1}  \cdots  \left( {\sum\limits_{{\ell_2} = 2}^{{\ell_3} - 1} {\left( {\sum\limits_{{\ell_1} = 1}^{{\ell_{2 - 1}}} \frac{1}
{
\left(\ell_1-\frac{1}{2}\right)^2\cdot
\cdots
\left(\ell_M-\frac{1}{2}\right)^2
} } \right)} } \right) \cdots } \right)}\right]x^{2M}
\]
$ \\ $
Die in den Summen auftretenden Brüche
\[
\frac{1}{\left(\ell_k-\frac{1}{2}\right)^2}\;;\;k\in\mathbb \{1,\ldots,n\}
\]
können wir wie folgt umformen
\[
\frac{1}{\left(\ell_k-\frac{1}{2}\right)^2}
=\frac{2^2}{\left(2\,\ell_k-1\right)^2}\;;\;k\in\mathbb \{1,\ldots,n\}
\]
Wir erhalten dann die Darstellung
\[
1-\left[\sum_{\ell_1=1}^{M}
\frac{2^2}{\left(2\,\ell_1-1\right)^2}
\right] x^2
\]
\[
+\left[\sum _{\ell_{{2}}=2}^{M}\left( \sum _{\ell_{{1}}=1}^{\ell_{{2}}-1}
\frac{2^2}{\left(2\,\ell_1-1\right)^2}
\cdot
\frac{2^2}{\left(2\,\ell_2-1\right)^2}
\right)\right] x^4
\]
\[
-\left[\sum _{\ell_{{3}}=3}^{M}
\left(\sum_{\ell_{{2}}=2}^{\ell_{{3}}-1} \left( 
\sum _{\ell_{{1}}=1}^{\ell_{{2}}-1}
\frac{2^2}{\left(2\,\ell_1-1\right)^2}
\cdot
\frac{2^2}{\left(2\,\ell_2-1\right)^2}
\cdot
\frac{2^2}{\left(2\,\ell_3-1\right)^2}
\right)\right)\right] x^6\;
\]
\[
+\left[\sum_{\ell_{{4}}=4}^{M}
\left(\sum_{\ell_{{3}}=3}^{\ell_{{4}}-1}
\left(\sum_{\ell_{{2}}=2}^{\ell_{{3}}-1}
\left(\sum_{\ell_{{1}}=1}^{\ell_{{2}}-1}
\frac{2^2}{\left(2\,\ell_1-1\right)^2}
\cdot
\frac{2^2}{\left(2\,\ell_2-1\right)^2}
\cdot
\frac{2^2}{\left(2\,\ell_3-1\right)^2}
\cdot
\frac{2^2}{\left(2\,\ell_4-1\right)^2}
\right)\right)\right)\right]x^8-+\cdots
\]
\[
+\left(-1\right)^M
\left[\sum\limits_{{\ell_M}=M}^M{
\left({\sum\limits_{{\ell_{M-1}}=M-1}^{{\ell_M}-1}\cdots
\left({\sum\limits_{{\ell_2}=2}^{{\ell_3}-1}{
\left({\sum\limits_{{\ell_1}=1}^{{\ell_{2-1}}}
\frac{2^2}{\left(2\,\ell_1-1\right)^2}
\cdots
\frac{2^2}{\left(2\,\ell_{M}-1\right)^2}
}\right)}}\right)\cdots}\right)}\right]x^{2\,M}
\]
$ \\ \\ $
Die $2$-er Potenzen
\[2^2\,;\quad 2^2\cdot 2^2=2^4\,;\quad 2^2\cdot 2^2\cdot 2^2=2^6\,;\quad\ldots\quad\,;\,\underbrace{2^2\cdot\ldots\cdot 2^2}_{M\times\, 2^2}=2^{2\,M}
\]
können wir ausklammern und vor die Summen bringen, so dass wir erhalten
\[
1-2\,^2\cdot\left[\sum_{\ell_1=1}^{M}
\frac{1}{\left(2\,\ell_1-1\right)^2}
\right] x^2
\]
\[
+2\,^4
\left[\sum _{\ell_{{2}}=2}^{M}\left( \sum _{\ell_{{1}}=1}^{\ell_{{2}}-1}
\frac{1}{\left(2\,\ell_1-1\right)^2\cdot\left(2\,\ell_2-1\right)^2}
\right)\right] x^4
\]
\[
-2\,^6
\left[\sum _{\ell_{{3}}=3}^{M}
\left(\sum_{\ell_{{2}}=2}^{\ell_{{3}}-1} \left( 
\sum _{\ell_{{1}}=1}^{\ell_{{2}}-1}
\frac{1}{\left(2\,\ell_1-1\right)^2
\cdot
\left(2\,\ell_2-1\right)^2
\cdot
\left(2\,\ell_3-1\right)^2}
\right)\right)\right] x^6\;
\]
\[
+2\,^8
\left[\sum_{\ell_{{4}}=4}^{M}
\left(\sum_{\ell_{{3}}=3}^{\ell_{{4}}-1}
\left(\sum_{\ell_{{2}}=2}^{\ell_{{3}}-1}
\left(\sum_{\ell_{{1}}=1}^{\ell_{{2}}-1}
\frac{1}{\left(2\,\ell_1-1\right)^2}
\cdot
\frac{1}{\left(2\,\ell_2-1\right)^2
\cdot
\left(2\,\ell_3-1\right)^2
\cdot
\left(2\,\ell_4-1\right)^2}
\right)\right)\right)\right]x^8-+\cdots
\]
\[
+\left(-1\right)^M\cdot2^{\,2\,M}
\left[\sum\limits_{{\ell_M}=M}^M{
\left({\sum\limits_{{\ell_{M-1}}=M-1}^{{\ell_M}-1}\cdots
\left({\sum\limits_{{\ell_2}=2}^{{\ell_3}-1}{
\left({\sum\limits_{{\ell_1}=1}^{{\ell_{2-1}}}
\frac{1}{\left(2\,\ell_1-1\right)^2
\cdots
\left(2\,\ell_{M}-1\right)^2}
}\right)}}\right)\cdots}\right)}\right]x^{2\,M}
\]
$ \\ $
Wir überlegen uns was geschieht, wenn wir $M \to \infty $ streben lassen.
$ \\ \\ $
Auf der linken Seite erhalten wir ein unendliches Produkt
\[
\underset{M\to \infty}{\mathop{\lim }}\;\;
\prod_{k=1}^{M}\left(1-\frac{x^2}{\left(k-\frac{1}{2}\right)^2}\right)
=\prod_{k=1}^{\infty}\left(1-\frac{x^2}{\left(k-\frac{1}{2}\right)^2}\right),
\]
das wir als die Produktdarstellung von $\cos\left(\pi\,x\right)$ identifizieren:
\[
\cos\left(\pi x\right)=\prod_{k=1}^{\infty}\left(1-\frac{x^2}{\left(k-\frac{1}{2}\right)^2}\right)
\]
Die Potenzreihe für $\cos\left(\pi\,x\right)$ steht in jeder Formelsammlung.
\[
\cos\left(\pi\,x\right)=1-\frac{\pi^2}{2!}x^2+\frac{\pi^4}{4!}x^4-\frac{\pi^6}{6!}x^6+-\cdots
\]
Auf der rechten Seite erhalten wir für $M \to \infty $ eine unendliche Potenzreihe, deren Koeffizienten wiederum aus unendlichen Summen bestehen.
\[
1-2^2\cdot\left[\sum_{\ell_1=1}^{\infty}
\frac{1}{\left(2\,\ell_1-1\right)^2}
\right] x^2
\]
\[
+2\,^4
\left[\sum _{\ell_{{2}}=2}^{\infty}\left( \sum _{\ell_{{1}}=1}^{\ell_{{2}}-1}
\frac{1}{\left(2\,\ell_1-1\right)^2\cdot\left(2\,\ell_2-1\right)^2}
\right)\right] x^4
\]
\[
-2\,^6
\left[\sum _{\ell_{{3}}=3}^{\infty}
\left(\sum_{\ell_{{2}}=2}^{\ell_{{3}}-1} \left( 
\sum _{\ell_{{1}}=1}^{\ell_{{2}}-1}
\frac{1}{\left(2\,\ell_1-1\right)^2
\cdot
\left(2\,\ell_2-1\right)^2
\cdot
\left(2\,\ell_3-1\right)^2}
\right)\right)\right] x^6\;
\]
\[
+2\,^8
\left[\sum_{\ell_{{4}}=4}^{\infty}
\left(\sum_{\ell_{{3}}=3}^{\ell_{{4}}-1}
\left(\sum_{\ell_{{2}}=2}^{\ell_{{3}}-1}
\left(\sum_{\ell_{{1}}=1}^{\ell_{{2}}-1}
\frac{1}{\left(2\,\ell_1-1\right)^2}
\cdot
\frac{1}{\left(2\,\ell_2-1\right)^2
\cdot
\left(2\,\ell_3-1\right)^2
\cdot
\left(2\,\ell_4-1\right)^2}
\right)\right)\right)\right]x^8\;-+\cdots
\]
\section*{3  Von der cosinus Produktdarstellung zu den Pi-Reihen}
Wir haben somit für $\cos\left(\pi\,x\right)$ zwei verschiedene Darstellungen mit Potenzreihen. Einmal die klassische
\[
\cos\left(\pi\,x\right)=1-\frac{\pi^2}{2!}x^2+\frac{\pi^4}{4!}x^4-\frac{\pi^6}{6!}x^6+-\cdots
\]
Zum anderen eine Darstellung durch ineinander geschachtelte Summen. Gemäß dem Eindeutigkeitssatz für Potenzreihen müssen beide identisch sein. Es muss also gelten:
\[
1-\frac{\pi^2}{2!}x^2+\frac{\pi^4}{4!}x^4-\frac{\pi^6}{6!}x^6
-\frac{\pi^8}{8\,!}x^8-+\cdots
\]
\[
=1-2\,^2\left[\sum_{\ell_1=1}^{\infty}
\frac{1}{\left(2\,\ell_1-1\right)^2}
\right] x^2
\]
\[
+2\,^4
\left[\sum _{\ell_{{2}}=2}^{\infty}\left( \sum _{\ell_{{1}}=1}^{\ell_{{2}}-1}
\frac{1}{\left(2\,\ell_1-1\right)^2\cdot\left(2\,\ell_2-1\right)^2}
\right)\right] x^4
\]
\[
-2\,^6
\left[\sum _{\ell_{{3}}=3}^{\infty}
\left(\sum_{\ell_{{2}}=2}^{\ell_{{3}}-1} \left( 
\sum _{\ell_{{1}}=1}^{\ell_{{2}}-1}
\frac{1}{\left(2\,\ell_1-1\right)^2
\cdot
\left(2\,\ell_2-1\right)^2
\cdot
\left(2\,\ell_3-1\right)^2}
\right)\right)\right] x^6\;
\]
\[
+2\,^8
\left[\sum_{\ell_{{4}}=4}^{\infty}
\left(\sum_{\ell_{{3}}=3}^{\ell_{{4}}-1}
\left(\sum_{\ell_{{2}}=2}^{\ell_{{3}}-1}
\left(\sum_{\ell_{{1}}=1}^{\ell_{{2}}-1}
\frac{1}{\left(2\,\ell_1-1\right)^2}
\cdot
\frac{1}{\left(2\,\ell_2-1\right)^2
\cdot
\left(2\,\ell_3-1\right)^2
\cdot
\left(2\,\ell_4-1\right)^2}
\right)\right)\right)\right]x^8\;-+\cdots
\]
$ \\ $
Falls zwei Potenzreihen identisch sind, müssen ihre Koeffizienten übereinstimmen.
$ \\ $
$ \\ $
Setzen wir die Koeffizienten von $x^2$ auf linker und rechter Seite gleich, so erhalten wir
\[
\frac{\pi^2}{2!}
=2^2\sum_{\ell_1=1}^{\infty}\frac{1}{\left(2\,\ell_1-1\right)^2}
\]
Nach Division durch $2^2$ ergibt sich
\[
\frac{\pi^2}{2\,!\cdot 2^2}
=\frac{1}{2\,!}\left(\frac{\pi}{2}\right)^2
=\frac{\pi^2}{8}
=\sum_{\ell_1=1}^{\infty}\frac{1}{\left(2\,\ell_1-1\right)^2}
=1+\frac{1}{3^2}+\frac{1}{5^2}+\frac{1}{7^2}+\frac{1}{9^2}+\ldots
\]
Diese Reihe taucht öfter im Zusammenhang mit Lösungen des Basler Problem auf.
$ \\ $
$ \\ $
Verfahren wir mit den Koeffizienten von $x^4$ genau so, so erhalten wir auf der linken Seite ein rationales Vielfaches von $\pi^4$ and auf der rechten Seite zwei ineinander geschachtelte Summen
\[
\frac{\pi^4}{4\,!}
=2\,^4\sum_{\ell_{{2}}=2}^{\infty}\left(\sum _{\ell_{{1}}=1}^{\ell_{{2}}-1}
\frac{1}{\left(2\,\ell_1-1\right)^2}\cdot
\frac{1}{\left(2\,\ell_2-1\right)^2}
\right)
\]
Das lässt sich wieder vereinfachen zu
\[
\frac{\pi^4}{4\,!\,2^4}
=\frac{1}{4\,!}\left(\frac{\pi}{2}\right)^4
=\sum_{\ell_{{2}}=2}^{\infty}\left(\sum _{\ell_{{1}}=1}^{\ell_{{2}}-1}\frac{1}{
\left(2\,\ell_1-1\right)^2\cdot
\left(2\,\ell_2-1\right)^2}
\right)
\]
Hierbei wollen wir unter dem Grenzwert der Summe
\[
\sum_{\ell_{{2}}=2}^{\infty}
\left(\sum _{\ell_{{1}}=1}^{\ell_{{2}}-1}\frac{1}{
\left(2\,\ell_1-1\right)^2\cdot
\left(2\,\ell_2-1\right)^2}
\right)
\]
$ \\ $
den Grenzwert der Partialsummen verstehen
\[
\sum_{\ell_2=2}^{2}
\left(\sum_{\ell_1=1}^{\ell_2-1}
\frac{1}
{
\left(2\ell_1-1\right)^2
\left(2\ell_2-1\right)^2
}
\right)
=\frac{1}{1^2\cdot 3^2}
\]
\[
\sum_{\ell_2=2}^{3}
\left(\sum_{\ell_1=1}^{\ell_2-1}
\frac{1}
{
\left(2\ell_1-1\right)^2
\left(2\ell_2-1\right)^2
}
\right)
=\frac{1}{1^2\cdot 3^2}+\frac{1}{1^2\cdot 5^2}+\frac{1}{3^2\cdot 5^2}
\]
\[
\sum_{\ell_2=2}^{4}
\left(\sum_{\ell_1=1}^{\ell_2-1}
\frac{1}
{
\left(2\ell_1-1\right)^2
\left(2\ell_2-1\right)^2
}
\right)
=\frac{1}{1^2\cdot 3^2}+\frac{1}{1^2\cdot 5^2}+\frac{1}{1^2\cdot7^2}
+\frac{1}{3^2\cdot 5^2}+\frac{1}{3^2\cdot7^2}
+\frac{1}{5^2\cdot7^2}
\]
\centerline{\vdots}
$ \\ $
Für die Koeffizienten von $x^6$ erhalten wir ebenfalls ein  rationales Vielfaches von  $\pi^6$, dargestellt durch eine dreifach verschachtelte Summe.
\[
\frac{\pi^6}{6\,!}
=2\,^6\frac{1}{6\,!}\left(\frac{\pi}{2}\right)^6
=\sum_{\ell_{{3}}=3}^{\infty}
\left(\sum_{\ell_{{2}}=2}^{\ell_{{3}}-1}
\left(\sum_{\ell_{{1}}=1}^{\ell_{{2}}-1}
\frac{1}{
\left(2\,\ell_1-1\right)^2\cdot
\left(2\,\ell_2-1\right)^2\cdot
\left(2\,\ell_3-1\right)^2}
\right)\right)
\]
Auch das können wir vereinfachen
\[
\frac{\pi^6}{6\,!\cdot 2^6}
=\frac{1}{6\,!}\left(\frac{\pi}{2}\right)^6
=\sum_{\ell_{{3}}=3}^{\infty}
\left(\sum_{\ell_{{2}}=2}^{\ell_{{3}}-1}
\left(\sum_{\ell_{{1}}=1}^{\ell_{{2}}-1}
\frac{1}{
\left(2\,\ell_1-1\right)^2\cdot
\left(2\,\ell_2-1\right)^2\cdot
\left(2\,\ell_3-1\right)^2}
\right)\right)
\]
Unter der Konvergenz verstehen wir wieder die Konvergenz der Partialsummen
\[
\sum_{\ell_3=3}^{3}
\left(\sum_{\ell_2=2}^{\ell_3-1}
\left(\sum_{\ell_1=1}^{\ell_2-1}
\frac{1}
{
\left(2\ell_1-1\right)^2
\left(2\ell_2-1\right)^2
\left(2\ell_3-1\right)^2
}
\right)\right)
=\frac{1}{1^{2} \cdot 3^{2}\cdot 5^{2}}
\]
\centerline{}
\[
\sum_{\ell_3=3}^{4}
\left(\sum_{\ell_2=2}^{\ell_3-1}
\left(\sum_{\ell_1=1}^{\ell_2-1}
\frac{1}
{
\left(2\ell_1-1\right)^2
\left(2\ell_2-1\right)^2
\left(2\ell_3-1\right)^2
}
\right)\right)
\]
\[
=\frac{1}{1^{2} \cdot 3^{2} \cdot 5^{2}}+\frac{1}{1^{2} \cdot 3^{2} \cdot 7^{2}}+\frac{1}{1^{2} \cdot 5^{2} \cdot 7^{2}}+\frac{1}{3^{2} \cdot 5^{2} \cdot 7^{2}}
\]
\centerline{}
\[
\sum_{\ell_3=3}^{5}
\left(\sum_{\ell_2=2}^{\ell_3-1}
\left(\sum_{\ell_1=1}^{\ell_2-1}
\frac{1}
{
\left(2\ell_1-1\right)^2
\left(2\ell_2-1\right)^2
\left(2\ell_3-1\right)^2
}
\right)\right)
\]
\[
=\frac{1}{1^{2} \cdot 3^{2} \cdot 5^{2}}+\frac{1}{1^{2} \cdot 3^{2} \cdot 7^{2}}+\frac{1}{1^{2} \cdot 3^{2} \cdot 9^{2}}+\frac{1}{1^{2} \cdot 5^{2} \cdot 7^{2}}+\frac{1}{1^{2} \cdot 5^{2} \cdot 9^{2}}
\]
\[
+\frac{1}{1^{2} \cdot 7^{2} \cdot 9^{2}}+\frac{1}{3^{2} \cdot 5^{2} \cdot 7^{2}}+\frac{1}{3^{2} \cdot 5^{2} \cdot 9^{2}}+\frac{1}{3^{2} \cdot 7^{2} \cdot 9^{2}}+\frac{1}{5^{2} \cdot 7^{2} \cdot 9^{2}}
\]
\centerline{\vdots}
$ \\ $
Als weitere Reihen erhalten wir, ohne jetzt auf Details einzugehen
\[
\frac{1}{8!}\left(\frac{\pi}{2}\right)^8=
\]
\[=
\sum_{\ell_4=4}^{\infty}
\left(\sum_{\ell_3=3}^{\ell_4-1}
\left(\sum_{\ell_2=2}^{\ell_3-1}
\left(\sum_{\ell_1=1}^{\ell_2-1}
\frac{1}
{
\left(2\ell_1-1\right)^2
\left(2\ell_2-1\right)^2
\left(2\ell_3-1\right)^2
\left(2\ell_4-1\right)^2
}
\right)\right)\right)
\]
\centerline{}
\[
\frac{1}{10!}\left(\frac{\pi}{2}\right)^{10}=
\]
\[
=\sum_{\ell_5=5}^{\infty}
\left(\sum_{\ell_4=4}^{\ell_5-1}
\left(\sum_{\ell_3=3}^{\ell_4-1}
\left(\sum_{\ell_2=2}^{\ell_3-1}
\left(\sum_{\ell_1=1}^{\ell_2-1}
\frac{1}
{
\left(2\ell_1-1\right)^2
\left(2\ell_2-1\right)^2
\left(2\ell_3-1\right)^2
\left(2\ell_4-1\right)^2
\left(2\ell_5-1\right)^2
}
\right)\right)\right)\right)
\]
\centerline{\vdots}
$ \\ $
Unter Konvergenz wollen wir stets die Konvergenz der Partialsummen verstehen.
$ \\ $
Die Gleichsetzung von Koeffizienten, wie wir es bisher für einzelne Potenzen getan haben, können wir auch für beliebige geradzahlige Exponenten vornehmen. Wir erhalten dann
\[
\frac{\pi^{2\,n}}{\left(2\,n\right)!}=
\]
\[
=2\,^{2\,n}\sum\limits_{{\ell_n}= n}^\infty {\left( {\sum\limits_{{\ell_{n-1}}=n-1}^{{\ell_n}-1} \cdots \left({\sum\limits_{{\ell_2}=2}^{{\ell_3}-1} {\left( {\sum\limits_{{\ell_1}=1}^{{\ell_{2-1}}}\frac{1}
{
\left(2\ell_1-1\right)^2
\left(2\ell_2-1\right)^2
\cdots
\left(2\ell_{n-1}-1\right)^2
\left(2\ell_n-1\right)^2
}}\right)}}\right)\cdots}\right)}
\]
$ \\ $
Wir dividieren die Gleichung noch durch $2\,^{2\,n}$ und erhalten die folgende Darstellung
\[
\frac{\pi^{2\,n}}{\left(2\,n\right)!}\cdot\frac{1}{2\,^{2\,n}}=
\frac{1}{\left(2\,n\right)!}\left(\frac{\pi}{2}\right)^{2\,n}=
\]
\[=\sum\limits_{{\ell_n}= n}^\infty {\left( {\sum\limits_{{\ell_{n-1}}=n-1}^{{\ell_n}-1} \cdots \left({\sum\limits_{{\ell_2}=2}^{{\ell_3}-1} {\left( {\sum\limits_{{\ell_1}=1}^{{\ell_{2-1}}}\frac{1}
{
\left(2\ell_1-1\right)^2
\left(2\ell_2-1\right)^2
\cdots
\left(2\ell_{n-1}-1\right)^2
\left(2\ell_n-1\right)^2
}}\right)}}\right)\cdots}\right)}
\]
$ \\ $
Schließlich müssen wir uns noch mit der Konvergenz der Reihen befassen. Zunächst stellen wir fest, dass nur positive Summanden auftreten. Somit ist das Majorantenkriterium die erste Wahl. Wir müssen also Vergleichsreihen finden, deren Summanden stets größer oder höchsten gleich der vorgelegten Reihen sind und die konvergieren.
$ \\ $
$ \\ $
Die Konvergenz der Reihe
\[
\sum _{\ell_{{1}}=1}^{\infty}\frac{1}{\left(2 l_{1}-1\right)^{2}}=\frac{\pi^2}{8}
\]
ist wohl bekannt. Es ist nichts weiter zu tun.
$ \\ $
$ \\ $
Für höhere $n$ verwenden wir Vergleichsreihen, deren Konvergenz bekannt ist \citep{KK},\, \cite{RC}
\[
\sum _{\ell_{{2}}=1}^{\infty}
\left(\sum_{\ell_{{1}}=1}^{\infty} 
\frac{1}
{
{\left(2\ell_1-1\right)}^2\cdot
{\left(2\ell_2-1\right)}^2
}\right)
=\left(\frac{\pi^2}{8}\right)^2;
\]
\[
\sum _{\ell_{{3}}=1}^{\infty}
\left(\sum_{\ell_{{2}}=1}^{\infty}
\left(\sum_{\ell_{{1}}=1}^{\infty}
\frac{1}{
{\left(2\ell_1-1\right)}^2\cdot
{\left(2\ell_2-1\right)}^2\cdot
{\left(2\ell_3-1\right)}^2
}\right)\right)
=\left(\frac{\pi^2}{8}\right)^3;
\]
\[
\sum_{\ell_{{4}}=1}^{\infty}
\left(\sum_{\ell_{{3}}=1}^{\infty}
\left(\sum_{\ell_{{2}}=1}^{\infty}
\left(\sum_{\ell_{{1}}=1}^{\infty}
\frac{1}{
\left(2\ell_1-1\right)^2\cdot
\left(2\ell_2-1\right)^2\cdot
\left(2\ell_3-1\right)^2\cdot
\left(2\ell_4-1\right)^2}
\right)\right)\right)
=\left(\frac{\pi^2}{8}\right)^4;
\]
\[
\vdots
\]
\[
\sum\limits_{{\ell_n} = 1}^{\infty} {\left( {\sum\limits_{{\ell_{n - 1}} =1}^{\infty}  \cdots  \left( {\sum\limits_{{\ell_2} = 1}^{\infty} {\left( {\sum\limits_{{\ell_1} = 1}^{\infty} {\frac{1}{
{\left(2\ell_1-1\right)^2 \cdots
\left(2\ell_n-1\right)^2}}} } \right)} } \right) \cdots } \right)}
=\left(\frac{\pi^2}{8}\right)^n\]
$ \\ \\ $
Wir vergleichen sie mit den hergeleiteten Reihen und stellen folgendes fest
\[
\sum _{\ell_{{2}}=2}^{\infty} \left( \sum _{\ell_{{1}}=1}^{\ell_{{2}}-1} 
\frac{1}{\left(2\ell_1-1\right)^2 \cdot
\left(2\ell_2-1\right)^2}\right)
<\sum _{\ell_{{2}}=1}^{\infty} \left( \sum _{\ell_{{1}}=1}^{\infty} 
\frac{1}{\left(2\ell_1-1\right)^2 \cdot
\left(2\ell_2-1\right)^2}\right)=\left(\frac{\pi^2}{8}\right)^2;
\]
\begin{align*}
&\sum _{\ell_{{3}}=3}^{\infty} \left( \sum _{\ell_{{2}}=2}^{\ell_{{3}}-1} \left( \sum _
{\ell_{{1}}=1}^{\ell_{{2}-1}} \frac{1}{
\left(2\ell_1-1\right)^2 \cdot
\left(2\ell_2-1\right)^2 \cdot
\left(2\ell_3-1\right)^2}\right)\right)\\ 
<
&\sum _{\ell_{{3}}=1}^{\infty} \left( \sum _{\ell_{{2}}=1}^{\infty} \left( \sum _
{\ell_{{1}}=1}^{\infty} \frac{1}{
\left(2\ell_1-1\right)^2\cdot
\left(2\ell_2-1\right)^2\cdot
\left(2\ell_3-1\right)^2}
\right)\right)=\left(\frac{\pi^2}{8}\right)^3;
\end{align*}
\centerline{}
\begin{align*}
&\sum_{\ell_{{4}}=4}^{\infty}
\left(\sum_{\ell_{{3}}=2}^{\ell_4-1}
\left(\sum_{\ell_{{2}}=2}^{\ell_3-1}
\left(\sum_{\ell_{{1}}=1}^{\ell_2-1}
\frac{1}{
\left(2\ell_1-1\right)^2\cdot
\left(2\ell_2-1\right)^2\cdot
\left(2\ell_3-1\right)^2\cdot
\left(2\ell_4-1\right)^2}
\right)\right)\right)
\\
<
&\sum_{\ell_{{4}}=1}^{\infty}
\left(\sum_{\ell_{{3}}=1}^{\infty}
\left(\sum_{\ell_{{2}}=1}^{\infty}
\left(\sum_{\ell_{{1}}=1}^{\infty}
\frac{1}{
\left(2\ell_1-1\right)^2\cdot
\left(2\ell_2-1\right)^2\cdot
\left(2\ell_3-1\right)^2\cdot
\left(2\ell_4-1\right)^2}
\right)\right)\right)
=\left(\frac{\pi^2}{8}\right)^4;
\end{align*}
\[
\vdots
\]
Für beliebiges $n\in\mathbb N$ erhalten wir
\begin{align*}
&\sum\limits_{{\ell_n} = n}^{\infty} {\left( {\sum\limits_{{\ell_{n - 1}} = n - 1}^{{\ell_n} - 1}  \cdots  \left( {\sum\limits_{{\ell_2} = 2}^{{\ell_3} - 1} {\left( {\sum\limits_{{\ell_1} = 1}^{{\ell_{2 - 1}}} {\frac{1}{{
\left(2\ell_1-1\right)^2\cdots
\left(2\ell_n-1\right)^2}}}}\right)}}\right)\cdots}\right)}
\\
<
&\sum\limits_{{\ell_n} = 1}^{\infty} {\left( {\sum\limits_{{\ell_{n - 1}} =1}^{\infty}  \cdots  \left( {\sum\limits_{{\ell_2} = 1}^{\infty} {\left( {\sum\limits_{{\ell_1} = 1}^{\infty} {\frac{1}{{
\left(2\ell_1-1\right)^2\cdots
\left(2\ell_n-1\right)^2
}}} } \right)} } \right) \cdots } \right)}
=\left(\frac{\pi^2}{8}\right)^n
\end{align*}
Somit sind die zur Disposition stehenden Reihen durch $\left(\frac{\pi^2}{8}\right)^n $ beschränkt und hieraus folgt nach dem Majorantenkriterium unmittelbar ihre Konvergenz.

\end{document}